\documentstyle[12pt]{article}
\def\beq{\begin{equation}}
\def\eeq{\end{equation}}
\def\bea{\begin{eqnarray}}
\def\eea{\end{eqnarray}}
\def\nn{\nonumber}
\def\ba{\begin{array}}
\def\ea{\end{array}}
\def\one{1\hskip -1mm{\rm l}}
\font\twelvemsbm=msbm10 at 12 true pt
\font\eightmsbm=msbm8
\font\sevenmsbm=msbm7
\newfam\msbmfam
\textfont\msbmfam=\twelvemsbm
\scriptfont\msbmfam=\eightmsbm
\scriptscriptfont\msbmfam=\sevenmsbm
\def\Bb#1{{\fam\msbmfam\relax#1}} 

\def\Z{\Bb Z}

\begin{document}
\baselineskip14pt

\begin{flushright}
q-alg/9602020
\end{flushright}

\noindent
{\large \bf The dual $(p,q)$-Alexander-Conway Hopf algebras and  
the associated universal ${\cal T}$-matrix}

\medskip

\noindent
R. Chakrabarti$^{\dag}$ and 
R. Jagannathan$^{\ddag}$\footnote{E-mail: jagan@imsc.ernet.in}

\vspace{.25cm}

\noindent
$^{\dag}$ {\small Department of Theoretical Physics, University of Madras, 
Guindy Campus, Madras-600025, {\phantom{IN\,}}INDIA} \\
$^{\ddag}$ {\small The Institute of Mathematical Sciences, C.I.T. Campus, 
Tharamani, Madras-600113, INDIA}

\vspace{1.5cm}

\baselineskip18pt

\noindent
{\bf Abstract}: The dually conjugate Hopf algebras $Fun_{p,q}(R)$ and  
$U_{p,q}(R)$ associated with the two-parametric $(p,q)$-Alexander-Conway 
solution $(R)$ of the Yang-Baxter equation are studied.  Using the Hopf 
duality construction, the full Hopf structure of the quasitriangular 
enveloping algebra $U_{p,q}(R)$ is extracted.  The universal ${\cal T}$-matrix 
for $Fun_{p,q}(R)$ is derived.  While expressing an arbitrary group element of 
the quantum group characterized by the noncommuting parameters in a 
representation independent way, the ${\cal T}$-matrix generalizes the familiar 
exponential relation between a Lie group and its Lie algebra.  The universal 
${\cal R}$-matrix and the FRT matrix generators, $L^{(\pm )}$, for $U_{p,q}(R)$ 
are derived from the ${\cal T}$-matrix.   

\vspace{1cm}

\noindent
{\bf 1. Introduction}
                                 
\renewcommand{\theequation}{1.{\arabic{equation}}}
\setcounter{equation}{0}

\medskip

\noindent
The quantum Yang-Baxter equation (QYBE) admits nonstandard 
solutions~\cite{LCS}-\cite{GSX} characterizing quasitriangular 
Hopf algebras which are not deformations of classical algebras.  The 
nonstandard quantum algebras associated with the Alexander-Conway 
solution of the QYBE has been studied~(\cite{LCS}-\cite{GSX}) and 
its two-parametric generalization has been obtained~\cite{BBCH,CJ1}.  
These algebras have the interesting property that, by the technique 
of superization~\cite{MR}, they can be associated with the graded 
quantized universal enveloping algebra (QUEA) 
$U_{\langle q \rangle}(gl(1|1))$, where $\langle q \rangle$ represents 
the set of pertinent deformation parameters.  Interestingly, the quantum 
group realization of the Alexander polynomial was obtained~\cite{KS,RS} 
from the algebra $U_q(gl(1|1))$ while yielding a free fermion model for 
the invariant.  In addition, for the nongeneric values of the 
deformation parameters, the nonstandard  $R$-matrices may be engendered 
using a coloured, generalized boson realization~\cite{GSX,CJ1} of the 
universal ${\cal R}$-matrices of the corresponding standard ungraded QUEA  
$U_{\langle q \rangle}(gl(2))$.  

The FRT-construction~\cite{FRT} associates to any solution $R$ of 
the QYBE a quantum matrix pseudogroup defined by the transfer matrix $T$.  
The elements of $T$ generate the function algebra 
$Fun_{\langle q \rangle}(R)$.  The QUEA $U_{\langle q \rangle}(R)$, dually 
conjugate to $Fun_{\langle q \rangle}(R)$ in the Hopf sense, may now be 
obtained~\cite{MR} if $Fun_{\langle q \rangle}(R)$ contains a group-like, 
`quantum determinent' type, element and some suitable ansatz for the 
matrix generators $L^{\pm}$ of $U_{\langle q \rangle}(R)$ exists.  The 
significance of a key notion capping the Hopf duality structure~\cite{FRT} 
was recently highlighted by Fr{\o}nsdal and Galindo~\cite{FG} in the context 
of $Fun_{p,q}(GL(2))$ and $U_{p,q}(gl(2))$~:  In a representation independent 
way, they derived a closed expression for the dual form 
\beq %1.1 
{\cal T} = \sum\,e^A E_A = {\cal T}_{e,E}\,, 
\label{calT}
\eeq
called the universal ${\cal T}$-matrix, and established it as the quantum group 
generalization of the familiar exponential map obtaining in the case of 
classical groups (see also~\cite{FG2,F}).  Here, the sets $\left\{ e^A \right\}$ 
and $\left\{ E_A \right\}$ are the respective basis elements of the dual Hopf 
algebras ${\cal A} = Fun_{\langle q \rangle}(R)$ and ${\cal U} = 
U_{\langle q \rangle}(R)$, satisfying the relation  
\beq %1.2 
\left\langle e^A\,,\,E_B \right\rangle = \delta^A_B\,, 
\label{eE}
\eeq
where $\langle\,,\,\rangle$ is a doubly nondegenerate bilinear form.   
Conversely, as in the classical case, an infinitesimal analysis of the quantum 
group elements would lead to the quantum algebra, as has been demonstrated by 
Finkelstein~\cite{Fi} in the example of $Fun_q(GL(2))$.  When the algebras 
$Fun_{\langle q \rangle}(R)$ and $U_{\langle q \rangle}(R)$ are finitely 
generated, a closed  expression of the ${\cal T}$-matrix may be obtained in 
terms of the two sets of generators.  The ${\cal T}$-matrix expresses a 
representation-free realization of a quantum group element depending on the 
noncommuting group parameters.  As pointed out in~(\cite{FG}-\cite{F}), the 
main usefulness of the ${\cal T}$-matrix derive from the fact that the 
transition matrices of the integrable models appear, upon specialization, in 
passing from the structure to the representations.  

The technique adopted in~(\cite{FG}-\cite{F}) may be summarized as follows.  The 
structure relations embodying the duality between the conjugate Hopf algebras 
may be expressed succinctly in terms of the ${\cal T}$-matrix as 
\bea %1.3
{\cal T}_{e,E} {\cal T}_{e',E} & = & {\cal T}_{\Delta (e),E}\,, \quad 
{\cal T}_{e,E} {\cal T}_{e,E'} = {\cal T}_{e,\Delta (E)}\,, 
\nn \\ 
{\cal T}_{\epsilon (e),E} & = & \one \,, \quad 
{\cal T}_{e,\epsilon (E)} = \one \,,
\nn \\
{\cal T}_{S(e),E} & = & {\cal T}^{-1}\,, \quad 
{\cal T}_{e,S(E)} = {\cal T}^{-1}\,, 
\label{TeE}
\eea
where $e$ and $e'$ ($E$ and $E'$) refer to two identical copies of 
$Fun_{\langle q \rangle}(R)$ ($U_{\langle q \rangle}(R)$) and $\Delta$, 
$\epsilon$ and $S$ denote, respectively, the coproduct, counit and the 
antipode maps.  Using the known Hopf structure of $Fun_{\langle q \rangle}(R)$, 
the dual Hopf structure of the QUEA $U_{\langle q \rangle}(R)$ may be read 
from~(\ref{TeE}).  When both the Hopf algebras are finitely generated, the 
basis elements $\left\{ e^A \right\}$ and $\left\{ E_A \right\}$ may be 
expressed as ordered monomials in the respective sets of generators.  A closed 
expression of the dual form may now be explicitly obtained in terms of the two 
sets of generators.  

Using the $q$-exponentials, an explicit expression for the ${\cal T}$-matrix 
corresponding to the dual Hopf algebras $Fun_{p,q}(GL(2))$ and $U_{p,q}(gl(2))$ 
was first obtained in~\cite{FG}.  The dual forms for the standard quantum 
$gl(n)$ and the twisted quantum $gl(n)$ were considered in~\cite{FG2} 
and~\cite{F}, respectively.  Following the Fr{\o}nsdal-Galindo approach, 
Bonechi {\em et al.}~\cite{Bo-et} derived the ${\cal T}$-matrices for some 
inhomogeneous quantum groups.  Morozov and Vinet~\cite{MV} constructed these 
generalized exponential maps for all standard simple quantum groups with a 
single deformation parameter.  The dual super-Hopf algebras 
$Fun_{p,q}(GL(1|1))$ and $U_{p,q}(gl(1|1))$ were studied in~\cite{CJ2} using 
the above approach.  In another developement, the $(p,q)$-generalization of the 
Wigner $d$-functions have been obtained~\cite{JV} using the finite dimensional 
representations of the ${\cal T}$-matrix for the algebra $Fun_{p,q}(GL(2))$.  
With a view to provide further concrete examples of the universal 
${\cal T}$-matrix we consider here the Hopf duality structure of the 
$(p,q)$-Alexander-Conway algebras and obtain the associated universal 
${\cal T}$-matrix.   

The universal ${\cal R}$-matrix for the $(p,q)$-Alexander-Conway algebra is 
known~\cite{CJ1} and we show here that it may be obtained from the universal 
${\cal T}$-matrix via a homomorphic map $\phi~: Fun_{p,q}(R) \longrightarrow 
U_{p,q}(R)$ using the relation
\beq %1.4 
( {\rm id} \otimes \phi ) {\cal T} = {\cal R}\,.
\label{TR}
\eeq
Such a map was first used in~\cite{F} in the context of quantum $gl(n)$.  
It is then obvious that in view of the link between the universal 
${\cal R}$-matrix and the FRT matrix generators $L^{(\pm )}$~\cite{FRT} it 
should be possible to obtain $L^{(\pm )}$ directly from the ${\cal T}$-matrix.    
Such a procedure has already been demonstrated in~\cite{VJ} in the case of 
$U_q(sl(2))$.  We shall exhibit here a similar derivation of 
$L^{(\pm )}$-matrices for $U_{p,q}(R)$ from the ${\cal T}$-matrix.  

\bigskip

\noindent
{\bf 2. Hopf structure of $Fun_{p,q}(R)$} 

\renewcommand{\theequation}{2.{\arabic{equation}}}
\setcounter{equation}{0}

\medskip

\noindent
We study the Hopf algebra associated with the two-parametric 
nonstandard solution~\cite{BBCH,CJ1} of the QYBE, namely, 
\beq %2.1
R = \left( 
\ba{cccc}
Q & 0 & 0 & 0 \\
0 & \lambda^{-1} & s & 0 \\
0 & 0 & \lambda & 0 \\
0 & 0 & 0 & -Q^{-1}  
\ea \right)\,,
\qquad s = Q-Q^{-1}\,.
\label{R}
\eeq
The defining relation of the quantum inverse scattering method~\cite{FRT},    
\beq %2.2
\sum_{m,n = 1,2}\,R_{im,jn} T_{mk} T_{nl} = \sum_{m,n = 1,2}\,T_{jn} T_{im} 
R_{mk,nl}\,, 
\label{RTT}
\eeq
with the $R$-matrix given by~(\ref{R}), describes a transfer matrix 
\beq %2.3 
T = \left[ t_{ij} \right] = \left(
\ba{cc}
a & b \\
c & d 
\ea \right)\,,
\label{T}
\eeq 
whose variable elements obey the braiding relations 
\bea %2.4 
ab & = & p^{-1}ba\,, \quad ac = q^{-1}ca\,, \quad 
db = -p^{-1}bd\,, \quad dc = -q^{-1}cd\,, \nn \\
p^{-1}bc & = & q^{-1}cb\,,  \quad ad-da = (p^{-1}-q)bc\,, 
\quad b^2 = 0\,, \quad c^2 = 0\,,
\label{comrel}
\eea
where
\beq %2.5  
p = \lambda Q\,, \quad \quad q = \lambda^{-1}Q\,.
\label{qlambda}
\eeq

If the diagonal elements of $T$, $a$ and $d$, are invertible, the 
elements $\{ a,b,c,d,$ $a^{-1},$ $d^{-1} \}$ generate a Hopf algebra 
$Fun_{p,q}(R)$ whose coalgebraic structure readily follows.  The 
coproduct, counit and the antipode maps 
are, respectively, given by 
\bea %2.6-8  
\Delta (T) & = & T \dot{\otimes} T\,, \nn \\
\Delta (a^{-1}) & = & a^{-1} \otimes a^{-1} - 
a^{-1}ba^{-1} \otimes a^{-1}ca^{-1}\,, \nn \\
\Delta (d^{-1}) & = & d^{-1} \otimes d^{-1} - 
d^{-1}cd^{-1} \otimes d^{-1}bd^{-1}\,, 
\label{delta} \\ 
\epsilon (T) & = & \one 
\label{epsilon} \\ 
S(T) & = & T^{-1}\,, \quad 
S(a^{-1}) = a - bd^{-1}c\,, \quad S(d^{-1}) = d - ca^{-1}b\,, 
\label{s}
\eea  
where
\beq %2.9  
T^{-1} = \left(
\ba{cc}
a^{-1} + a^{-1}bd^{-1}ca^{-1} & -a^{-1}bd^{-1} \\
-d^{-1}ca^{-1} & d^{-1} + d^{-1}ca^{-1}bd^{-1} 
\ea \right)
\label{t^-1}
\eeq 
and $\dot{\otimes}$ denotes the tensor product coupled with matrix 
multiplication.  Despite the appearance of the relations 
$(b^2 = 0$,  $c^2 = 0)$ in~(\ref{comrel}), suggestive of a 
superalgebraic structure, the Hopf algebra $Fun_{p,q}(R)$ is bosonic as 
it follows the tensor product rule 
\beq %2.10
(\Gamma_1 \otimes \Gamma_2)(\Gamma_3 \otimes \Gamma_4) = 
\Gamma_1\Gamma_3 \otimes \Gamma_2\Gamma_4\,, \quad 
\forall \ \ \Gamma \in Fun_{p,q}(R)\,. 
\label{tensor}
\eeq

In the Hopf algebra $Fun_{p,q}(R)$ an invertible group-like element $D$ 
exists: 
\beq %2.11
D = ad^{-1} - bd^{-1}cd^{-1}\,, \qquad 
D^{-1} = da^{-1} - ba^{-1}ca^{-1}\,.
\label{D}
\eeq
Using~(\ref{comrel}), the commutation relations for $D$ follow:
\beq %2.12
[D,a] = 0\,, \quad [D,d\,] = 0\,, 
\quad \{ D,b \} = 0\,, \quad \{ D,c \} = 0\,.
\label{dcomrel}
\eeq
The induced coalgebra maps for $D$ are 
\beq %2.13
\Delta (D) = D \otimes D\,, \quad 
\epsilon (D) = 1\,, \quad 
S(D) = D^{-1}\,,
\label{dcoalg}
\eeq
as obtained from the relations~(\ref{delta}-\ref{t^-1}). 

A Gauss decomposition of the $T$-matrix~(\ref{T}) as 
\beq %2.14
T = \left(
\ba{cc}
1 & 0 \\
\zeta & 1 
\ea \right) \left( 
\ba{cc}
a & 0 \\
0 & \hat{d} 
\ea \right) \left(  
\ba{cc}
1 & \xi \\
0 & 1 
\ea \right)
\label{factor}
\eeq
introduces new variables related to the old ones by 
\beq %2.15
b = a\xi\,, \qquad c = \zeta a\,, \qquad d = \zeta a \xi + \hat{d}\,,
\label{newvar}
\eeq
where $\hat{d}$ is an invertible element.  The element $D$ and its 
inverse now read 
\beq %2.16
D = a\hat{d}^{-1}\,, \quad \quad D^{-1} = \hat{d}a^{-1}\,.
\label{DD^-1}
\eeq
The algebra~(\ref{comrel}) assumes the form 
\bea %2.17
a\xi & = & p^{-1}\xi a\,, \quad a\zeta = q^{-1}\zeta a\,, \quad 
\hat{d}\xi = -p^{-1}\xi \hat{d}\,, \quad 
\hat{d}\zeta = -q^{-1}\zeta \hat{d}\,, \nn \\ {}
[a,\hat{d}\,] & = & 0\,, \quad [\xi ,\zeta ] = 0\,, \quad
\xi^2 = 0\,, \quad \zeta^2 = 0\,.
\label{newalg}
\eea 
The coalgebra maps~(\ref{delta}-\ref{t^-1}) are rewritten as 
\bea %2.18-20
\Delta (a) & = & a \otimes a + a\xi \otimes \zeta a\,, \qquad 
\Delta (\xi ) = \one \otimes \xi + \xi \otimes D^{-1}\,, \nn \\
\Delta (\zeta ) & = & \zeta \otimes \one + D^{-1} \otimes \zeta\,, 
\qquad \Delta (\hat{d}) = \hat{d} \otimes \hat{d} 
- \hat{d}\xi \otimes \zeta \hat{d}\,,
\label{newdelta} \\ 
\epsilon (a) & = & 1\,, \quad \epsilon (\xi ) = 0\,, \quad 
\epsilon (\zeta ) = 0\,, \quad \epsilon (\hat{d}) = 1\,,
\label{newepsilon} \\
S(a) & = & a^{-1} + \xi a^{-1}D\zeta\,, \qquad S(\xi ) = -\xi D\,, 
\nn \\ 
S(\zeta ) & = & -D\zeta\,, \qquad 
S(\hat{d}) = \hat{d}^{-1} - \xi \hat{d}^{-1}D\zeta\,.
\label{news}
\eea

Assuming that the algebra can be augmented with the logarithms of $a$ and 
$\hat{d}$, we use the map    
\beq %2.21
a = {\rm e}^x\,, \quad \quad \hat{d} = {\rm e}^{\hat{x}}\,,
\label{xhat}
\eeq
and a reparametrization 
\beq %2.22
p = {\rm e}^{-\omega}\,, \quad \quad q = {\rm e}^{-\nu} 
\label{omeganu}
\eeq
to convert~(\ref{newalg}) to the algebraic structure   
\bea %2.23
[x,\xi ] & = & \omega \xi\,, \quad [x,\zeta ] = \nu \zeta\,, \quad 
[\hat{x},\xi ] = \hat{\omega}\xi\,, \quad 
[\hat{x},\zeta ] = \hat{\nu}\zeta\,, \nn \\  
  &   & \qquad \qquad \qquad \hat{\omega} = \omega +{\rm i}\pi\,, \quad     
\hat{\nu} = \nu +{\rm i}\pi \nn \\ {}  
[x, \hat{x} ] & = & 0\,, \quad [\xi ,\zeta ] = 0\,, \quad
\xi^2 = 0\,, \quad \zeta^2 = 0\,,
\label{lie}
\eea
where the choice of identical phases for $\hat{\omega}$ and 
$\hat{\nu}$ is necessitated by the requirement that the single 
deformation parameter limit $(p$ $=$ $q)$ exists and, via a duality 
construction, yields the $q$-Alexander-Conway algebra~(\cite{LCS}-\cite{GSX}).  
Following the technique of Ref.~\cite{FG}, the coalgebra maps for~(\ref{lie}) 
are seen to be  
\bea %2.24-26
\Delta (x) & = & x \otimes \one + \one \otimes x + 
{{\nu +\omega}\over{{\rm e}^{\nu}-{\rm e}^{-\omega}}} 
\xi \otimes \zeta\,, \nn \\
\Delta (\xi ) & = & \one \otimes \xi + 
\xi \otimes {\rm e}^{\hat{x}-x}\,, \nn \\
\Delta (\zeta ) & = & \zeta \otimes \one 
+ {\rm e}^{\hat{x}-x} \otimes \zeta\,,
\nn \\
\Delta (\hat{x}) & = & \hat{x} \otimes \one + \one \otimes \hat{x} +
{{\nu + \omega +{\rm i}2\pi}\over{{\rm e}^{\nu}-{\rm e}^{-\omega}}} 
\xi \otimes \zeta\,,
\label{deltalie} \\ 
\epsilon (x) & = & \epsilon (\xi ) = \epsilon (\zeta ) 
= \epsilon (\hat{x}) = 0\,,
\label{epsilonlie} \\ 
S(x) & = & -x + {{\nu + \omega}\over{{\rm e}^{\nu}-{\rm e}^{-\omega}}}
\xi {\rm e}^{x-\hat{x}}\zeta\,, \qquad 
S(\xi ) = -\xi {\rm e}^{x-\hat{x}}\,, \nn \\
S(\zeta ) & = & -{\rm e}^{x-\hat{x}}\zeta\,, \qquad 
S(\hat{x}) = -\hat{x} + {{\nu + \omega + {\rm i}2\pi}\over
{{\rm e}^{\nu}-{\rm e}^{-\omega}}}\xi {\rm e}^{x-\hat{x}}\zeta\,.
\label{slie}
\eea 
The above Hopf structure suggests that $Fun_{p,q}(R)$ may be embedded 
in the enveloping algebra of $(x, \hat{x}, \xi , \zeta )$ 
satisfying the algebraic relations~(\ref{lie}) and endowed with a 
noncocommutative coproduct rule.  This enveloping algebra is dual to 
the Hopf algebra $U_{p,q}(R)$.  

All the Hopf algebra axioms for the coalgebra maps~(\ref{deltalie}-\ref{slie}) 
can be explicitly proved.  To this end, the following identities are to be 
noted.  If $X$ and $\tilde{X}$ are elements satisfying 
\beq %2.27
[x,\xi ] = c_{\xi} \xi \,, \quad 
[X,\zeta ] = c_{\zeta} \zeta\,, \quad 
[X,\tilde{X} ] = 0 
\label{XX}
\eeq
and $C$ is an arbitrary $c$-number, we have 
\bea %2.28
\exp ( X \otimes \one + \one \otimes X + C \xi \otimes \zeta )
  & = & \left( {\rm e}^X \otimes \one \right) \left( \one \otimes \one 
+ \tilde{C} \xi \otimes \zeta \right) \left( \one \otimes {\rm e}^X \right)\,, 
\label{eX} \\ 
\exp \left( -X + \xi \tilde{X} \zeta \right) 
  & = & {\rm e}^X + \Lambda \xi {\rm e}^{-X} \tilde{X} \zeta 
\label{e-x}
\eea
where $\Lambda = \frac{{\rm e}^{c_{\zeta}} - {\rm e}^{-c_{\xi}}}{c_{\zeta} + 
c_{\xi}}$, $\tilde{C} = C \Lambda$.  From the above identities and~(\ref{lie}), 
it follows that the one parametric set of elements 
\beq %2.30
D(\theta ) = {\rm e}^{\theta (x - \hat{x} )} 
\label{dtheta}
\eeq
satisfies a group-like property 
\beq %2.31
\Delta (D(\theta )) = D(\theta ) \otimes D(\theta )\,, \qquad 
S (D(\theta )) = D(\theta )^{-1} 
\label{group}
\eeq
for any $\theta \in \Z$ .  It is interesting to note that the elements 
$\{ D(\theta ) | \theta \in \Z \}$ form a discrete group that cannot be 
embedded in a continuous one-parameter group.    

Before proceeding to derive the structure of $U_{p,q}(R)$ and construct the 
universal ${\cal T}$-matrix, let us note an interesting application of the 
characterization of the elements of the $T$-matrix~(\ref{T}) through the 
algebraic algebraic structure~(\ref{lie}).  Let 
$\left\{ \left( a_1\,,\,a_1^\dagger \right)\,,\,\left( a_2\,,\,a_2^\dagger 
\right) \right\}$ and $\left\{ \left( b_1\,,\,b_1^\dagger \right)\,,\,\left( 
b_2\,,\,b_2^\dagger \right) \right\}$ be fermion operators corresponding to 
two disparate, commuting, Fermi fields; they may be the basic operators from 
which a pair of $2$-nd order para-Fermi operators can be obtained through 
the well known Green's ansatz~\cite{G}, or they can be considered as 
belonging to two Fermi fields of different colours.  Such commuting Fermi 
fields can also be constructed from a single Fermi field~\cite{JVa}.  The 
algebra of these $(a,b)$-operators is 
\bea %2.32
\left\{ a_i , a_j^\dagger \right\} & = & \delta_{ij}\,, \quad 
\left\{ a_i , a_j \right\} = 0\,, \quad i,j = 1,2\,, 
\nn \\
\left\{ b_i , b_j^\dagger \right\} & = & \delta_{ij}\,, \quad 
\left\{ b_i , b_j \right\} = 0\,, \quad i,j = 1,2\,, 
\nn \\
\left[ a_i , b_j \right] & = & 0\,, \quad 
\left[ a_i , b_j^\dagger \right] = 0\,, \quad i,j = 1,2\,.
\label{ab}
\eea
Using the Bogoliubov transformation, we define
\bea %2.33
\xi & = & ( \cos \alpha ) a_1 + ( \sin \alpha ) a_2^\dagger \quad 
\left(~{\rm or}~~( \cos \alpha ) a_2 - ( \sin \alpha ) a_1^\dagger~\right)\,,  
\nn \\ 
\zeta & = & ( \cos \beta ) b_1 + ( \sin \beta ) b_2^\dagger \quad 
\left(~{\rm or}~~( \cos \beta ) b_2 - ( \sin \beta ) b_1^\dagger~\right)\,.  
\label{xizeta}
\eea 
Then, we can take  
\beq %2.34
x = \omega \xi^\dagger \xi + \nu \zeta^\dagger\zeta + f(\alpha , \beta)\,, 
\qquad
\hat{x} = \hat{\omega} \xi^\dagger \xi + \hat{\nu} \zeta^\dagger\zeta + 
g(\alpha , \beta)\,, 
\label{xx}
\eeq 
where $f(\alpha , \beta)$ and $g(\alpha , \beta)$ are two arbitrary $c$-number 
functions of $\alpha$ and $\beta$.  Now, it is seen that, 
using~(\ref{newvar}),~(\ref{omeganu}) and~(\ref{ab})-(\ref{xx}), it is 
possible to realize the variable group element~(\ref{T}) parametrized by the 
classical variables $\alpha$ and $\beta$.  A similar realization of 
$GL_{p,q}(1|1)$ in terms of a fermion field was obtained in~\cite{CJ2} 
following~\cite{JV} in which $GL_{p,q}(2)$ was realized in terms of a boson 
field.   
 
\bigskip  

\noindent
{\bf 3. Dual Hopf algebra $U_{p,q}(R)$ and the universal 
${\cal T}$-matrix}

\renewcommand{\theequation}{3.{\arabic{equation}}}
\setcounter{equation}{0}

\medskip

\noindent
The monomials $\{ e^A|e^A = \zeta^{a_1}x^{a_2}\hat{x}^{a_3}\xi^{a_4}$,  
$A = (a_1,a_2,a_3,a_4)$, $a_1,a_4 = (0,1)$, $a_2,a_3 = 0,1,2,\ldots \}$ 
constitute a basis for ${\cal A} = Fun_{p,q}(R)$ obeying the multiplication and 
the induced coalgebra maps  
\bea %3.1,2
e^Ae^B & = & \sum_C\,f^{AB}_C e^C\,,
\label{meAeB} \\ 
\Delta (e^A) & = & \sum_{BC}\,h^A_{BC} e^B \otimes e^C\,, \quad 
\epsilon \left( e^A \right) = \delta^A_{\underline{0}}\,, \quad 
S\left( e^A \right) = \sum_B S^A_B e^B\,.
\label{coeA}
\eea
The unit element is obtained by choosing $A = \underline{0}$, where 
$\underline{0} = (0,0,0,0)$.  The elements $\{ E_A \}$, defined by~(\ref{eE}), 
form a basis set for the algebra ${\cal U}$, dual to ${\cal A}$, namely, the 
QUEA $U_{p,q}(R)$.  The duality construction~(\ref{TeE}) enforces the 
following Hopf structures for the basis set $\{ E_A \}$~:  %3.3,4
\bea %3.3-4 
E_AE_B & = & \sum_C\,h^C_{AB} E_C\,,
\label{mEAEB} \\  
\Delta (E_A) & = & \sum_{BC}\,f^{BC}_A E_B \otimes E_C\,.
\label{deltaEA} 
\eea

Using the algebra~(\ref{lie}) the structure tensor $f^{AB}_C$ is 
derived:
\bea %3.5
f^{AB}_C & = & 
\bar{\delta}^{a_1b_1}\bar{\delta}^{a_4b_4}\delta^{a_1+b_1}_{c_1}
\theta^{a_2+b_2}_{c_2}\theta^{a_3+b_3}_{c_3}\delta^{a_4+b_4}_{c_4} 
\nn \\
  &   &  \ \ \ \times \sum_{kl}\,
\left(
\ba{c}
a_2 \\
k
\ea \right)\left(
\ba{c}
b_2 \\
c_2-k
\ea \right)\left(
\ba{c}
a_3 \\
l
\ea \right)\left(
\ba{c}
b_3 \\
c_3-l
\ea \right) \nn \\
  &   &  \ \ \ \times (\nu b_1)^{a_2-k}(-\omega a_4)^{b_2-c_2+k}
((\nu + {\rm i}\pi )b_1)^{a_3-l} 
(-(\omega + {\rm i}\pi )a_4)^{b_3-c_3+l}\,,
\label{f}
\eea 
where $\bar{\delta}^{ab} = 
\delta^a_0\delta^b_0+\delta^a_1\delta^b_0+\delta^a_0\delta^b_1$ and 
$\theta^a_b = 1\,(0)$ if $a \geq b\ (\leq b)$.  The tensor $h^A_{BC}$ is 
determined using the induced coproduct for the basis set $\{e^A\}$, namely,  
\beq %3.6 
\Delta (e^A) = \Delta (\zeta )^{a_1} \Delta (x)^{a_2}
\Delta (\hat{x})^{a_3} \Delta (\xi)^{a_4}\,.
\label{delta-eA}
\eeq 
The following special cases, necessary for determining the dual
algebraic structure, may be directly read from~(\ref{delta-eA}):
\renewcommand{\theequation}{3.7{\alph{equation}}}
\setcounter{equation}{0}
\bea %3.7a-h
h^A_{B \underline{0}} & = & \delta^A_B\,, \quad \qquad
h^A_{\underline{0} B} = \delta^A_B\,,
\label{3.7a} \\ 
h^A_{0 b_2 b_3 00 c_2 c_3 0} & = & \delta^{a_1}_0 
\delta^{a_2}_{b_2 + c_2} \delta^{a_3}_{b_3 + c_3} \delta^{a_4}_0
\left(
\begin{array}{c}
a_2 \\
b_2 \\
\end{array}
\right)
\left(
\begin{array}{c}
a_3 \\
b_3 \\
\end{array}
\right)\,,
\label{3.7b} \\ 
h^A_{1000B}  & = & 
\delta^{a_1}_1 \delta^0_{b_1} 
\prod_{i=2}^{4} \delta^{a_i}_{b_i}\,, \qquad
h^A_{B0001}  = 
\left( \prod_{i=1}^{3} {\delta}^{a_i}_{b_i} \right)
\delta^{a_4}_1 \delta^0_{b_4}\,,
\label{3.7c} \\ 
h^A_{01001000} & = & 
\delta^{a_1}_1 \delta^{a_2}_1 \delta^{a_3}_0 \delta^{a_4}_0 -
\delta^{a_1}_1 \prod_{i=2}^{4} \delta^{a_i}_0\,,
\label{3.7d} \\ 
h^A_{00101000} & = & 
\delta^{a_1}_1 \delta^{a_2}_0 \delta^{a_3}_1 \delta^{a_4}_0 + 
\delta^{a_1}_1 \prod_{i=2}^{4} \delta^{a_i}_0\,,
\label{3.7e} \\ 
h^A_{00010100} & = &  
\delta^{a_1}_0 \delta^{a_2}_1 \delta^{a_3}_0 \delta^{a_4}_1 -
\left( \prod_{i=1}^{3} \delta^{a_i}_0 \right) \delta^{a_4}_1\,, 
\label{3.7f} \\ 
h^A_{00010010} & = &  
\delta^{a_1}_0 \delta^{a_2}_0 \delta^{a_3}_1 \delta^{a_4}_1 +
\left( \prod_{i=1}^{3} \delta^{a_i}_0 \right) \delta^{a_4}_1\,,
\label{3.7g} \\ 
h^A_{00011000} & = &  
\delta^{a_1}_1 \delta^{a_2}_0 \delta^{a_3}_0 \delta^{a_4}_1 +
\delta^{a_1}_0 \delta^{a_4}_0 \Omega_{a_2,a_3}\,,
\label{3.7h}
\eea 
where
\renewcommand{\theequation}{3.{\arabic{equation}}}
\setcounter{equation}{7}
\beq %3.8 
\Omega_{j,k} =
\frac{ \nu^j ( \nu + {\rm i} \pi )^k - ( - \omega )^j(- \omega -{\rm i} 
\pi)^k} { {\rm e}^\nu - {\rm e}^{- \omega}}\,.
\label{Omega}
\eeq
The antipode maps for the basis elements $\{e^A\}$ are obtained 
from~(\ref{slie}) ~:
\bea %3.9
S(e^A) 
& = & S(\xi)^{a_4} 
S( \hat{x} )^{a_3}
S(x)^{a_2}
S(\zeta )^{a_1} 
\nn \\ 
& = &
(-1)^{a_2+a_3} \zeta^{a_1}
\left( \hat{x} + a_1 (\nu + {\rm i} \pi) - a_4 (\omega + {\rm i} \pi) 
- \frac{\nu + \omega + {\rm i} 2 \pi}{ {\rm e}^\nu - {\rm e}^{-\omega}}
\zeta {\rm e}^{x - \hat{x}} \xi \right)^{a_3} \nonumber \\
&   & \ \ \ \times \left( x + a_1 \nu  - a_4 \omega
- \frac{\nu + \omega }{ {\rm e}^\nu - {\rm e}^{- \omega}} \zeta
{\rm e}^{x - \hat{x}} \xi \right)^{a_2}
{\rm e}^{(a_1 + a_4) ( x - \hat{x})} \xi^{a_4}\,.
\label{S-eA}
\eea 
The second equality in~(\ref{S-eA}) follows by using the commutation  
relations~(\ref{lie}) and will be later used to compute the antipode 
maps for the dual basis elements.

Employing the duality property, we now extract the multiplication
relations for the dual basis elements $\{E_A\}$. From~(\ref{mEAEB}) 
and~(\ref{3.7a}) the unit element follows directly:
\beq %3.10
E_A E_{\underline{0}} = E_A\,, \quad
E_{\underline{0}} E_A = E_A \ \ \longrightarrow
\ \ E_{\underline{0}} = 1\,.
\label{unit}
\eeq 
By choosing the generators of the dual algebra as
\beq %3.11
E_- = E_{1000}\,, \qquad
H = E_{0100}\,, \qquad
\hat{H} = E_{0010}\,, \qquad
E_+ = E_{0001}\,, \qquad
\label{generators}
\eeq 
and using the tensor structure~(\ref{3.7b}) and~(\ref{3.7c}) we express an 
arbitrary dual basis element as
\beq %3.12 
E_A = (a_2!a_3!)^{-1} E_-^{a_1} H^{a_2} \hat{H}^{a_3} E_+^{a_4}\,. 
\label{Es}
\eeq 
Further use of the special values of the structure tensor $h^A_{BC}$ 
in (3.7) now yields the commutation relations for the generators of the dual 
alglebra $U_{p,q}(R)$: 
\bea %3.13
[H, E_\pm ] & = & \pm E_\pm\,, \quad 
[\hat{H}, E_\pm ] = \mp E_\pm\,, \quad 
[H, \hat{H} ] = 0\,, \quad 
E_\pm ^2 = 0\,, \nn \\ {}  
[E_+ , E_- ] & = & \frac{ {\rm e}^{\nu (H+\hat{H}) } g^{-1} - 
{\rm e}^{-\omega (H+\hat{H})} g }{ {\rm e}^\nu - {\rm e}^{-\omega} }\,,
\label{upqr}
\eea 
where
\beq %3.14
g = {\rm e}^{-{\rm i} \pi \hat{H}}\,.
\label{g}
\eeq 
The algebra~(\ref{upqr}) may now be exploited to compute the general 
expression for the structure tensor $h^A_{BC}$, which reads
\bea %3.15 
h^A_{BC} & = & (-1)^{b_2+c_2-a_2} \bar{\delta}^{b_1c_1} 
\bar{\delta}^{b_4c_4} \delta^{b_1+c_1}_{a_1} \theta^{b_2+c_2}_{a_2}
\theta^{b_3+c_3}_{a_3}\delta^{b_4+c_4}_{a_4} \nn \\ 
  &   &  \ \times 
a_2!a_3!(b_2!b_3!c_2!c_3!)^{-1} \sum_{kl}
\left( \ba{cc}
b_2 \\
k
\ea \right) \left( \ba{cc}
c_2 \\
a_2-k 
\ea \right) \left( \ba{cc}
b_3 \\
l
\ea \right) \left( \ba{cc}
c_3 \\
a_3-l 
\ea \right) \nn \\
  &  & \ \times 
c_1^{b_2+b_3-k-l}b_4^{c_2+c_3-a_2-a_3+k+l} + 
\delta^{a_1}_{b_1} \delta^{b_4}_1 \delta^{c_1}_1 \delta^{a_4}_{c_4} 
a_2!a_3! \nn \\
  &   & \ \times
(b_2!c_2!b_3!c_3!(a_2-b_2-c_2)!(a_3-b_3-c_3)!)^{-1} 
\Omega_{a_2-b_2-c_2,a_3-b_3-c_3}\,.  
\label{genhabc}
\eea

The coproduct rules for the generators of the dual algebra $U_{p,q}(R)$ 
are obtained from~(\ref{deltaEA}) and~(\ref{f}): 
\bea %3.16
\Delta (H) & = & H \otimes \one + \one \otimes H\,, \qquad 
\Delta (\hat{H}) = \hat{H} \otimes \one + \one \otimes \hat{H}\,, 
\nn \\
\Delta (E_+) & = & E_+ \otimes {\rm e}^{-\omega(H+\hat{H})} g^{-1} + 
\one \otimes E_+\,, \quad  
\Delta (E_-) = E_- \otimes \one + {\rm e}^{\nu (H+\hat{H})} g \otimes E_-\,. 
\nn \\ 
   &     &    
\label{DeltaH}
\eea  
The counit maps for the dual generators follow from~(\ref{deltaEA})~:  
\beq %3.17
\epsilon (X) = 0\,, \qquad \forall\,X \in (H,\hat{H},E_\pm )\,.
\label{epsilonX}
\eeq 
Special values of the tensor $\left\{ S^A_B \right\}$, obtained 
from~(\ref{S-eA}), yields, via~(\ref{deltaEA}), the antipode maps for the 
dual generators~: 
\bea %3.18 
S(H) & = & -H\,, \qquad S(\hat{H}) = -\hat{H}\,, \nn \\
S(E_+) & = & g^{-1} {\rm e}^{\omega (H+\hat{H})} E_+\,,  \qquad 
S(E_-) = E_-{\rm e}^{-\nu (H+\hat{H})} g\,. 
\label{S(H)}
\eea

For the single deformation parameter case $(\nu$ $=$ $\omega)$, the 
Hopf structure~(\ref{upqr}, \ref{DeltaH}, \ref{epsilonX}, and 
\ref{S(H)}) for the dual algebra $U_{p,q}(R)$, after appropriate 
mappings, reduce to the results obtained in~\cite{MR} using the FRT 
construction~\cite{FRT}.  Following~\cite{MR}, we note that the element $g$ 
has a group-like coalgebra structure 
\beq %3.19
\Delta (g) = g \otimes g\,, \quad 
\epsilon (g) = \one \,, \quad
S(g) = g^{-1} 
\label{gco}
\eeq 
and the element $g^2$ is central.  By superization map~\cite{MR} the 
super-Hopf algebra $U_{p,q}(gl(1|1))$ may be realized from the quotient 
algebra $U_{p,q}(R)/g^2-\one$.  After a redefinition of the generators 
\beq %3.20 
Z = \frac{1}{2}(H+\hat{H})\,, \quad 
J = \frac{1}{2}(H-\hat{H})\,, \quad 
\chi_\pm = E_\pm Q^{\mp Z}\lambda^{-Z+\frac{1}{2}}\,, 
\eeq
\label{redefgen} 
the structure of the Hopf algebra $U_{p,q}(R)$ assumes the form 
\bea %3.21-24 
[J,\chi_\pm ] & = & \pm \chi_\pm\,, \quad 
[\chi_+ , \chi_-] = \frac{Q^{2Z}g-Q^{-2Z}g^{-1}}{Q-Q^{-1}}\,, \quad 
\chi_\pm ^2 = 0\,, \nn \\ {} 
[Z,X] & = & 0\,, \qquad \forall\,X \in (J,\chi_\pm )\,,
\label{newAlg} \\ 
\Delta (Z) & = & Z \otimes \one + \one \otimes Z\,, \qquad 
\Delta (J) = J \otimes \one + \one \otimes J\,, \nn \\
\Delta (\chi_+) & = & \chi_+ \otimes Q^Z \lambda^Z g + 
    Q^{-Z} \lambda^{-Z} \otimes \chi_+\,, \nn \\
\Delta (\chi_-) & = & \chi_- \otimes Q^Z \lambda^{-Z} + 
Q^{-Z} \lambda^Z g^{-1} \otimes \chi_-\,,
\label{newDelta} \\ 
\epsilon (X) & = & 0\,, \qquad \forall\,X \in (Z,J,\chi_\pm )\,,
\label{newEpsilon} \\ 
S(Z) & = & -Z\,, \quad S(J) = -J\,, \quad 
S(\chi_+) = g^{-1}\chi_+\,, \quad S(\chi_-) = \chi_- g\,.
\label{newS}
\eea 
In spite of the nonstandard commutation relations~(\ref{newAlg}) the 
Hopf algebra $U_{p,q}(R)$ is bosonic as it obeys the direct product 
rule~(\ref{tensor}).  The algebra $U_{p,q}(R)$ is quasitriangular with the 
universal ${\cal R}$-matrix given by~\cite{CJ1} 
\beq %3.25 
{\cal R} = (-1)^{(Z - J) \otimes (Z - J)} Q^{2(Z \otimes J 
+ J \otimes Z)} \lambda^{2(Z \otimes J - J \otimes Z)} 
{\rm e}^{s Q^Z \lambda^Z \chi_+ \otimes Q^{-Z} \lambda^Z \chi_-}\,.  
\label{calR}
\eeq 
For later use, we note that, in terms of the generators~(\ref{generators}),   
\beq %3.26 
{\cal R} =  (-1)^{\hat{H} \otimes \hat{H}} Q^{\left( H \otimes H - 
\hat{H} \otimes \hat{H} \right)} \lambda^{- \left( H \otimes \hat{H} - 
\hat{H} \otimes H \right)} 
{\rm e}^{\lambda s E_+ \otimes E_-} 
\label{calRgen}
\eeq

Finally, following the prescription~(\ref{calT}), we now explicitly write 
down the universal ${\cal T}$-matrix of $Fun_{p,q}(R)$: 
\beq %3.27 
{\cal T} =  {\rm e}^{\zeta E_-} {\rm e}^{x H + \hat{x} \hat{H}} 
{\rm e}^{\xi E_+}\,.  
\label{expliTau} 
\eeq 
It may be noted that in contrast to the case of $Fun_{p,q}(GL(2))$, discussed 
in~\cite{FG}, in the present case the $q$-exponentials do not appear in the 
expression for the ${\cal T}$-matrix as a consequence of the algebraic property 
$E_\pm^2 = 0$.  Corresponding to the two-dimensional irreducible 
representation of the generators of $U_{p,q}(R)$ given by 
\beq %3.28 
H = \left( \ba{cc}
1 & 0 \\
0 & 0 
\ea \right)\,, \quad 
\hat{H} = \left( \ba{cc}
0 & 0 \\
0 & 1 
\ea \right)\,, \quad 
E_+ = \left( \ba{cc}
0 & 1 \\
0 & 0 
\ea \right)\,, \quad 
E_- = \left( \ba{cc}
0 & 0 \\
1 & 0 
\ea \right)\,,
\label{2drep}
\eeq 
the universal ${\cal T}$-matrix~(\ref{expliTau}) reduces, as required, 
to the $T$-matrix~(\ref{T}). 

Let us now consider some applications of the ${\cal T}$-matrix~(\ref{expliTau}).  
First, we use the relation~(\ref{TR}) to reproduce the 
${\cal R}$-matrix~(\ref{calRgen}) from the ${\cal T}$-matrix~(\ref{expliTau}).  
To this end, the required homomorphism is  
\bea %3.29 
\phi (x) & = & (\ln Q) H - (\ln \lambda ) \hat{H}\,, \quad 
\phi (\hat{x}) = (\ln \lambda ) H - (\ln Q - {\rm i}\pi ) 
\hat{H}\,, 
\nn \\ 
\phi (\xi ) & = & \lambda s E_-\,, \quad 
\phi (\zeta ) = 0\,, 
\label{explihomo}
\eea 
which can be proved to satisfy the algebra~(\ref{lie}) using the  
relations~(\ref{upqr}).  The map~(\ref{TR}) now gives the 
${\cal R}$-matrix~(\ref{calRgen}) from the corresponding 
${\cal T}$-matrix~(\ref{expliTau}).  As noted by Fr{\o}nsdal~\cite{F}, in 
the context of quantum $gl(n)$, there should exist an alternative homomorphism, 
say $\tilde{\phi}~: {\cal A} \longrightarrow {\cal U}$, such that  
$({\rm id} \otimes \tilde{\phi}) {\cal T}$ corresponds to the alternative form 
of the universal ${\cal R}$-matrix.  We note that in the present case 
the alternative homomorphism is given by  
\bea %3.30 
\tilde{\phi} (x) & = & - (\ln Q) H - (\ln \lambda ) \hat{H}\,, 
\quad 
\tilde{\phi} (\hat{x}) = (\ln \lambda ) H + (\ln Q - {\rm i} \pi ) 
\hat{H}\,, 
\nn \\ 
\tilde{\phi} (\xi ) & = & 0\,, \quad 
\tilde{\phi} (\zeta ) = - \lambda s E_+\,.  
\label{homotilde}
\eea 
such that, as expected, we get 
\bea %3.31
( {\rm id} \otimes \tilde{\phi} ) {\cal T} & = & \tilde{\cal R} 
= {\rm e}^{-\lambda s E_- \otimes E_+} 
\lambda^{H \otimes \tilde{H} - \tilde{H} \otimes H} Q^{\tilde{H} \otimes 
\tilde{H} - H \otimes H} (-1)^{\tilde{H} \otimes \tilde{H}} 
\nn \\  
   & = & (\sigma ({\cal R}))^{-1} = \left( {\cal R}^{(+)} \right)^{-1}\,, 
\label{urtilde}
\eea 
where $\sigma (a \otimes b) = b \otimes a$. 

As is well known~\cite{FRT}, $L^{(\pm )}$, the FRT matrix generators of 
${\cal U}$, can be obtained from the universal ${\cal R}$-matrix.  Now, since 
the universal ${\cal R}$-matrix is related to the universal ${\cal T}$-matrix 
through a map it is obvious that we can get the $L^{(\pm )}$-matrices also 
directly from ${\cal T}$.  Let us demonstrate this procedure in the present 
case following~\cite{VJ} where the $L^{(\pm )}$-matrices of $U_q(sl(2))$ have 
been derived directly from the corresponding universal ${\cal T}$-matrix.     
Using the representation~(\ref{2drep}) in~(\ref{explihomo}) and 
(\ref{homotilde}) we get two representations of $\{ x, \hat{x}, \xi , \zeta 
\}$ as follows~:    
\bea %3.32-33 
\pi^+~: &   &  x = \left( 
\ba{cc}
\ln Q & 0 \\
0 & - \ln \lambda 
\ea \right)\,, \quad 
\hat{x} = \left( 
\ba{cc}
\ln \lambda & 0 \\
0 & - (\ln Q) + {\rm i}\pi  
\ea \right)\,, 
\nn \\
   &   & \xi = \left( 
\ba{cc}
0 & 0 \\
\lambda s & 0 
\ea \right)\,, \quad 
\zeta = 0\,, 
\\ 
\pi^-~: &   &  x = \left( 
\ba{cc}
- \ln Q & 0 \\
0 & - \ln \lambda 
\ea \right)\,, \quad 
\hat{x} = \left( 
\ba{cc}
\ln \lambda & 0 \\
0 & (\ln Q) - {\rm i}\pi  
\ea \right)\,, 
\nn \\
   &   & \xi = 0\,, \quad 
\zeta = \left( 
\ba{cc}
0 & - \lambda s \\
0 & 0 
\ea \right)\,.  
\label{2drepx}
\eea
Now, we obtain the required results~: 
\bea %3.34-35
(\pi^+ \otimes \one )({\cal T}) & = & L^{(+)} = \left( 
\ba{cc}
\lambda^{\hat{H}} Q^H & 0 \\
s g^{-1} \lambda^{1-H} Q^{-\hat{H}} E_+ & g^{-1} \lambda^{-H} Q^{-\hat{H}} 
\ea \right)\,, 
\\
(\pi^- \otimes \one )({\cal T}) & = & L^{(-)} = \left( 
\ba{cc}
\lambda^{\hat{H}} Q^{-H} & -s E_- \lambda^{1-H} Q^{\hat{H}} g \\
0 & \lambda^{-H} Q^{\hat{H}} g 
\ea \right)\,. 
\label{lpm}
\eea
It can be verified directly that these $L^{(\pm )}$-matrices generate the 
algebra $U_{p,q}(R)$~(\ref{upqr}) in the FRT-approach~\cite{FRT}.   

\bigskip

\noindent{\bf 4. Conclusion}

\medskip

We have studied the dually paired Hopf algebras 
${\cal A} = Fun_{p,q}(R)$ and ${\cal U} = U_{p,q}(R)$ associated with the 
nonstandard $R$-matrix~(\ref{R}) involving two independent parameters. Using 
the technique developed by Fr{\o}nsdal and Galindo~\cite{FG}-\cite{F} we have 
extracted the full Hopf structure of the algebra $U_{p,q}(R)$.  A 
representation of the quantum group element, known as the universal 
${\cal T}$-matrix, and associated with the corresponding dual form, is 
obtained and found to exhibit the suitably modified familiar exponential map 
relating the Lie group with the corresponding Lie algebra. 
 
Following the construction in~\cite{MV}, the present derivation of 
the ${\cal T}$-matrix may be extended to general nonstandard 
algebras~\cite{CL} related to the Hopf superalgebras 
$U_{\langle q \rangle}(gl(m|n))$.  The corresponding nonstandard $R$-matrices 
are known~\cite{VL} to underly the mathematical structure of the 
Park-Schultz (generalized six vertex) model and its associated 
quantum spin chains, which are of much current interest~\cite{MRi,SNW}.  

\bigskip

\noindent
{\bf Acknowledgement}

\medskip

\noindent
We are thankful to the referee for useful comments.

\end{document}